\documentclass[12pt, a4paper]{article}

\usepackage{a4wide}
\usepackage{amsmath}
\usepackage{amsfonts}
\usepackage{enumerate}
\date{}
\begin{document}
\title{\large\textbf{ON NONMEASURABLE SUBGROUPS IN\\ GENERALIZED MEASURABLE STRUCTURES }}
\author{\normalsize
\textbf{S.BASU \& D.SEN}\\
}
\date{}	
\maketitle
{\small \noindent \textbf{{\textbf{AMS subject classification ($2010$):}}}} {\small $28A05$, $28D05$, $03E05$, $03E50$.\\}
{\small\noindent \textbf{\textbf{\textmd{\textbf{{Keywords and phrases :}}}}} Continuum hypothesis, Banach-Kuratowski Matrix, admissible transfinite matrix, Polish groups, $\sigma$-algebra admissible with respect to a small system, small system upper semicontinuous relative to a $\sigma$-algebra.}\\
\vspace{.03cm}\\
{\normalsize
\textbf{\textbf{ABSTRACT:}} In this paper, we prove a result on nonmeasurable subgroups in commutative Polish groups with respect to more generalized structures than $\sigma$-finite measures.\\
\begin{center}
\section{\large{INTRODUCTION}}
\end{center}
\normalsize  It can be shown using Continuum hypothesis that there exists a countable family of subgroups of $\mathbb{R}$ such that no nonzero, $\sigma$-finite, diffused measure exists which includes all of these sets in its domain. Kharazishvili[6] extended these results in commutative groups. He proved the following two theorems :\\
\vspace{.01cm}

\textbf{THEOREM K1 : } Assume continuum hypothesis. Let $G$ be a commutative group having the cardinality of the continuum. Then there exists a countable family $(G_{i})_{_{i<\omega}}$ ($\omega$ is the first infinite cardinal) of subgroups of $G$, such that, for any non-zero, $\sigma$-finite, diffused measure $\mu$ on $G$, at least one of the groups is nonmeasurable with respect to $\mu$.\\
\vspace{.01cm}

\textbf{THEOREM K2 : } Let $G$ be any commutative group of cardinality $\omega_{1}$ (the first uncountable cardinal). Then there exists a countable family $(G_{i})_{_{i<\omega}}$ ($\omega$ is the first infinite cardinal) of subgroups of $G$, such that, for every non-zero, $\sigma$-finite, diffused measure $\mu$ on $G$, at least one of the groups $G_{i}$ is nonmeasurable with respect to $\mu$.\
\vspace{.05cm}

It may be noted that Theorem K2 does not require Continuum hypothesis.\
\vspace{.05cm}

In proving the above two theorems, Kharazishvili used two different type of matrices namely, the Banach-Kuratowski matrix [7] and the admissible transfinite matrix [7] over any set cardinality $\omega_{1}$.\\
\vspace{.01cm}

\textbf{DEFINITION 1.1 :} Given any set $E$ with card${(E)}= \omega_{_{1}}$, a double sequence $(E_{m}{_,}{_{n}})_{m<\omega}{_,}{_{n<\omega}}$ is called a Banach Kuratowski matrix if\\
(i) $E_{m}{_,}{_{0}}\subseteq E_{m}{_,}{_{1}}\subseteq \ldots\ldots\subseteq E_{m}{_,}{_{n}}\subseteq \ldots\ldots$ for every $m<\omega$.\\
(b) $E= \cup \{E_{m}{_,}{_{n}} : n< \omega\}$ and\\
(c) for any arbitrary chosen function $f:\omega\rightarrow \omega$, $E_{0}{_,}{_{f(0)}}\cap E_{1}{_,}{_{f(1)}}\cap \ldots\ldots E_{m}{_,}{_{f(m)}}\cap\ldots$ is atmost countable.\\
\vspace{.05cm}

In the family Let $F= \omega^{\omega}$ of all functions from $\omega$ into $\omega$, we define a preordering as follows:\\ $f\preceq g$ iff there exists a natural number $n=n(f,g)$ such that $f(m)\leq g(m)$ for all $m\geq n$. Then under the assumption of continuum hypothesis, a subset $E= \{f_{\xi} : \xi< \omega_{_{1}}\}$ can be defined which satisfies the following two conditions:\\
(a) The set $E$ is cofinal in $F$, i.e, for any arbitrary chosen member $f$ of $F$, there exists $\xi< \omega_{_{1}}$ such that $f\prec f_{\xi}$.\\
(b) For any two $\xi$ and $\rho$ related by the inequlity $\xi< \rho< \omega_{_{1}}$ the relation $f_{\rho}\preceq f_{\xi}$ does not hold.\
\vspace{.05cm}

Conditions (a) and (b) imply that card${(E)}= \omega_{_{1}}$ and if we put $E_{_{m,n}}= \{f_{\xi}: {f_{\xi}}(m)\leq {n}\}$, then the double sequence $(E_{m}{_,}{_{n}})_{m<\omega}{_,}{_{n<\omega}}$ forms a Banach-Kuratowski matrix over $E$.\
\vspace{.05cm}

It may be checked that there does not exist [7] any nonzero, $\sigma$-finite, diffused measure defined simultaneously for all the sets $E_{_{m,n}}$.\\
\vspace{.05cm}

\textbf{DEFINITION 1.2 :} Let card$(E)$= $\omega_{1}$. A double family $(E_{m}{_,}{_{\xi}})_{m<\omega}{_,}{_{\xi<\omega_{_{1}}}}$ of subsets of $E$ is called an admissible transfinite matrix for $E$ if\\
(i) For each $\xi<\omega_{1}$, the partial family $(E_{_{n,\xi}})_{_{n<\omega}}$ is increasing by inclusion and $\displaystyle{\bigcup_{n<\omega}}{E_{_{n,\xi}}}= E$\\
(ii) For any natural number $n$, there exists a natural number $m=m(n)$ such that for any set $\mathcal D\subseteq \omega _{1}$ having card$(\mathcal D)=m$, card$(\cap\{E_{_{n,\xi}}: \xi\in\mathcal D\})\leq\omega$.\
\vspace{.05cm}

For any set $E$ with card$(E)=\omega_{1}$, there always exists an admissible transfinite matrix for $E$. This follows from the existence of Ulam's matrix [7] for $E$. However, in the above two Definitions, we identify every infinite cardinal with the initial ordinal representing the cardinal.\
\vspace{.1cm}

In this paper, we obtain a modified generalization of the above two theorems in commutative Polish group (i.e in commutative groups whose underlying topology is separable and complete metrizable). Our development follows a pattern that is similar to that followed by Kharazishvili. The difference is that in our case the nonmeasurable subgroups are Bernstein. We also avoid using measures and instead replace it by a suitably designed measurable structure formed out of a $\sigma$-algebra which is admissible (in a sense defined in the following section) with respect to a small system of sets. It is worthwhile to mention here that the notion of small system was originally introduced by Ri$\acute{e}$can and Neubrunn [12] and susequently used by several authors in giving abstract formulations of many well-known results of classical measure and integration theory (see [4], [5], [10], [11]).\\
\vspace{.05cm}

\section{\large{PRELIMINARIES AND RESULTS}}
\normalsize In a topological space $X$, a set is called totally imperfect if it is contains no nonempty perfect set. A Bernstein set [7] is that which together with its complement are both totally imperfect. It was first constructed by Bernstein in 1908. Bernstein sets are nonmeasurable with respect to the completion of every nonzero $\sigma$-finite diffused Borel measure and this phenomenon also characterises a Bernstein set. Since in every uncountable Polish space the family of all nonempty perfect sets has the cardinality of the continuum, Bernstein sets can be constructed in such spaces.\
\vspace{.05cm}

Let $X$ be a nonempty set and $\mathcal S$ be a $\sigma$-algebra of subsets of $X$. $\mathcal S$ is called diffused if $\{x\}\in \mathcal S$ for every $x\in X$. By a small system on $X$ (the Definition is a modified version of the one given in [12]) we mean\\
\vspace{.01cm}

\textbf{DEFINITION 2.1: } A sequence $\{\mathcal N_{n}\}_{_{n<\omega}}$, where each $\mathcal N_{n}$ is a class of subsets of $X$ satisfying the following properties:\\
(i) $\emptyset\in \mathcal N_{n}$.\\
(ii) $E\in \mathcal N_{n}$ and $F\subseteq E$ implies $F\in \mathcal N_{n}$. In otherwords, each $\mathcal N_{n}$ is a hereditary class.\\
(iii) $E\in \mathcal N_{n}$ and $F\in{\displaystyle{\bigcap_{n=1}^{\infty}}\hspace{.01cm}{\mathcal N_{n}}}$ implies $E\cup F\in\mathcal N_{n}$.\\
(iv) For any $m= 1,2\ldots$, there exists a natural number $m^{\prime}> m$ such that for any one-to-one correspondence $k\rightarrow n_{_{k}}$ with $n_{k} > m^{\prime}$, ${\displaystyle{\bigcup_{k=1}^{\infty}{E_{n_{k}}}}}\in \mathcal N_{_{m}}$ whenever $E_{n_{k}}\in \mathcal N_{n_{k}}$.\\
(v)  For any $p,q$ there exists $m > p,q$  such that $\mathcal N_{_{m}} \subseteq \mathcal N_{ _{p}}, \mathcal N _{_{q}}$ In otherwords, the system is directed.\
\vspace{.03cm}

We further define a small system $\{\mathcal N_{n}\}_{_{n<\omega}}$ on $X$ to be diffused if $\{x\}\in \mathcal N_{n}$ for all $n$ and $x\in X$.\\
\vspace{.05cm}

\textbf{DEFINITION 2.2: } A $\sigma$-algebra $\mathcal S$ on $X$ is called admissible with respect to a small system $\{\mathcal N_{n}\}{_{_{n<\omega}}}$ if\\
(i) $\mathcal S\setminus\mathcal N_{n}\neq\emptyset\neq\mathcal S\cap\mathcal N_{n}$.\\
(ii) Each $\mathcal N_{n}$ has a $S$-base which means that each $E\in\mathcal N_{n}$ is contained in some $F\in {\mathcal S\cap\mathcal N_{n}}$.\\
(iii) $\mathcal S\setminus\mathcal N_{n}$ satisfies countable chain condition, which means that the cardinality of any arbitrary collection of mutually disjoint sets from $\mathcal S\setminus\mathcal N_{n}$ is atmost $\omega$.\
\vspace{.05cm}

More general Definitions than above are given in [1], [2] and [3]. We set $\mathcal N_{\infty}= {\displaystyle{\bigcap_{n<\omega}}}\hspace{.02cm}{\mathcal N_{n}}$. By virtue of conditions (ii), (iv) of Definition $2.1$, $\mathcal N_{\infty}$ is a $\sigma$-ideal. This $\sigma$-algebra $\mathcal S$ together with the $\sigma$-ideal $\mathcal N_{\infty}$ generates the $\sigma$-algebra $\mathcal{\widetilde S}$ whose members are of the form $(X\setminus Y)\cup Z$ where $X\in\mathcal S$ and $Y,Z\in \mathcal N_{\infty}$ and this gives rise to the measurable structure $(\mathcal{\widetilde S}, \mathcal N_{\infty})$. From admissibility of $\mathcal S$ with respect to $\{\mathcal N_{n}\}{_{_{n<\omega}}}$ it follows that $(\mathcal{\widetilde S}, \mathcal N_{\infty})$ satisfies countable chain condition.\\
\vspace{.05cm}

\textbf{DEFINITION 2.3 [4], [12]: } A small system $\{\mathcal N_{n}\}_{n<\omega}$ is said to be upper semicontinuous relative to a $\sigma$-algebra $\mathcal S$ if for every nested sequence $\{E_{n}\}_{n<\omega}$ of sets from $\mathcal S$ satisfying $E_{n}\notin \mathcal N_{m}$ for some $m$ and $n=1,2\ldots$,  ${\displaystyle{\bigcap_{n<\omega}}{E_{n}}}\notin \mathcal N_{\infty}$.\
\vspace{.05cm}

More general Definition than above is given in [3].\
\vspace{.05cm}

Below are some propositions which we need for our purpose.\\
\vspace{.05cm}

\textbf{PROPOSITION 2.4: } If $\{\mathcal N_{n}\}_{n<\omega}$ is upper semicontinuous relative to $\mathcal S$, then $\{\mathcal N_{n}\}_{n<\omega}$ is also upper semicontinuous relative to $\mathcal{\widetilde S}$.\\
\vspace{.05cm}

\textbf{PROOF : } Let $\{E_{n}\}_{n<\omega}$ be a nested sequence of sets from $\mathcal {\widetilde S}$ such that $E_{n}\notin \mathcal N_{m}$ for some $m$ and $n=1,2\ldots$. We write $E_{n}= F_{n}\Delta P_{n}$ where $F_{n}\in\mathcal S$ and $P_{n}\in \mathcal N_{\infty}$. As $\mathcal S$ is admissible with respect to $\{\mathcal N_{n}\}_{n<\omega}$, so $P_{n}\subseteq Q_{n}\in{\mathcal S\cap\mathcal N_{\infty}}$. Hence $\{F_{n}\setminus{\displaystyle{\bigcup_{k=1}^{n}{Q_{k}}}}\}_{n<\omega}$ forms a nested sequence with $F_{n}\setminus{\displaystyle{\bigcup_{k=1}^{n}{Q_{k}}}}\notin \mathcal N_{m}$ for $n=1,2\ldots$ by conditions (ii) and (iii) of Definition 2.1. But $F_{n}\setminus{\displaystyle{\bigcup_{i=1}^{n}{Q_{i}}}}\in\mathcal S$ and $\{\mathcal N_{n}\}_{n<\omega}$ is upper semicontinuous relative to $\mathcal S$, so $\displaystyle{\bigcap_{n<\omega}}(F_{n}\setminus{\displaystyle{\bigcup_{i=1}^{n}{Q_{i}}}})\notin \mathcal N_{\infty}$ which again implies that ${\displaystyle{\bigcap_{n<\omega}}{E_{n}}}\notin \mathcal N_{\infty}$. This proves the proposition.\\
\vspace{.05cm}

\textbf{PROPOSITION 2.5: } If $Y\notin \mathcal N_{\infty}$, then there exists a natural number $m$ such that no subset $M$ of $Y$ can belong to $\mathcal N_{m}$ if its complement in $Y$ i.e $Y\setminus M$ belong to $\mathcal N_{m}$.\\
\vspace{.05cm}

\textbf{PROOF :} Choose a natural number $r$ such that $Y\notin \mathcal N_{r}$. Since there exist $p,q$ (by condition (iv) of Definition 2.1) such that $\mathcal N_{p}\cup \mathcal N_{q}\subseteq \mathcal N_{r}$ and also $m$ (by condition (v) of Definition 2.1) such that $\mathcal N_{m}\subseteq \mathcal N_{p}$ and $\mathcal N_{m}\subseteq\mathcal N_{q}$, so if $M$ and $Y\setminus M$ both belong to $\mathcal N_{m}$, then their union belong to $\mathcal N_{p}\cup\mathcal N_{q}$ and so in $\mathcal N_{r}$- a contradiction.\\
\vspace{.05cm}

\textbf{PROPOSITION 2.6: } If a $\sigma$-algebra $\mathcal S$ is admissible with respect to a small system $\{\mathcal N_{n}\}_{n<\omega}$, then $\mathcal{\widetilde S}$ is also admissible with respect to $\{\mathcal N_{n}\}_{n<\omega}$.\\
\vspace{.05cm}

\textbf{PROOF :} The first two conditions of Definition 2.2 are obvious. Let $E= F\Delta P\in \mathcal {\widetilde S}\setminus \mathcal N_{n}$ where $F\in\mathcal S$ and $P\in\mathcal N_{\infty}$. By condition (ii) of Definition 2.2, there exists $Q\in\mathcal N_{\infty}\cap\mathcal S$ such that $P\subseteq Q$. Clearly $F\setminus Q\in\mathcal S\setminus\mathcal N_{n}$ for otherwise by condition (ii) and (iii) of Definition 2.1, $E\subseteq(F\setminus Q)\cup Q\in\mathcal N_{n}$. Thus every set in $\mathcal{\widetilde S}\setminus \mathcal N_{n}$ contains a set in $\mathcal S\setminus \mathcal N_{n}$. But $\mathcal S$ being admissible with respect to $\{\mathcal N_{n}\}_{n<\omega}$, $\mathcal S\setminus \mathcal N_{n}$ satisfies countable chain condition. Hence $\mathcal{\widetilde S}\setminus\mathcal N_{\infty}$ also satisfies countable chain condition which proves condition (iii) of Definition 2.2.\\
\vspace{.05cm}

\textbf{PROPOSITION 2.7: } Let $\{Z_{\alpha} : \alpha<\omega_{1}\}$ be an uncountable collection of sets from $\mathcal S$ and $m$ is a positive integer such that $\displaystyle{\bigcap_{\alpha\in \mathcal D}{Z_{\alpha}}}\in \mathcal N_{\infty}$ for every $m$-element subset $\mathcal D$ where $\mathcal S$ is admissible with respect to $\{\mathcal N_{n}\}_{n<\omega}$. Then there exists an uncountable subset $A$ of $\omega_{1}$ such that $Z_{\alpha}\in\mathcal N_{\infty}$ for every $\alpha\in A$.\
\vspace{.05cm}

A proof of the above proposition follows from the notion of admissibility and the inductive argument used in proving Lemma 4 (Chapter 13, [7]).\\
\vspace{.05cm}

\textbf{PROPOSITION 2.8: } Assuming Continuum hypothesis, any uncountable Polish group can be expressed as the direct sum of two Bernstein subgroups; i.e, if $(G,+)$ is an uncountable Polish group then $G= H+B$ where $H\cap B= \{0\}$ and $H,B$ are Bernstein subgroups of $G$.\\
\vspace{.05cm}

\textbf{PROOF :}  Since $G$ is an uncountable Polish group, card($G$)=$c$ (the cardinality of the continuum), which is again the cardinality of the class of all nonempty perfect subsets of $G$. By the virtue of continuum hypothesis, we may now consider an enumeration of this class $\{P_{\xi}: \xi<\omega_{1}\}$ and using transfinite recursion construct two $\omega_{1}$-sequences $\{x_{\xi}: \xi<\omega_{1}\}$ and $\{y_{\xi} : \xi<\omega_{1}\}$ such that (i) $x_{\xi}, y_{\xi}\in P_{\xi}$ for $\xi<\omega_{1}$ \\
(ii) the family $\{x_{\xi}: \xi<\omega_{1}\}\cup \{y_{\xi}: \xi<\omega_{1}\}$ is independent in the sense that if $m_{1}z_{i_{_{1}}}+ m_{2}z_{i_{_{2}}}+\ldots m_{k}z_{i_{_{k}}}= 0$ where $\{i_{1}, i_{2}\ldots i_{k}\}\subseteq \omega$; $m_{1}, m_{2}\ldots m_{k}$ are all integers and $\{z_{i_{1}}, z_{i_{2}}\ldots z_{i_{k}}\}\subseteq \{x_{\xi}: \xi<\omega_{1}\}\cup \{y_{\xi}: \xi<\omega_{1}\}$ then $m_{1}= m_{2}=\ldots = m_{k}= 0$.\
\vspace{.05cm}

Let $H$ be the subgroup generated by $\{x_{\xi}: \xi<\omega_{1}\}$ and $B$ be the maximal subgroup containing $\{y_{\xi}: \xi<\omega_{1}\}$ such that $H\cap B= \{0\}$. Then $H$ and $B$ are required Bernstein subgroups of $G$.\
\vspace{.05cm}

The above proof is similar to the proof of Lemma 1 (chapter 18, [8]).\
\vspace{.05cm}

Consider the ideal $\mathcal I$ generated by the family $\{\displaystyle\bigcup_{g\in\Sigma}g+B : \Sigma\hspace{.1cm} is\hspace{.1cm} a\hspace{.05cm} countable\hspace{.05cm} subset\hspace{.05cm} of\hspace{.05cm} H\}$. It is a $\sigma$-ideal, and also proper because for each member $Z$ of this family, there exists an uncountable family $\{g_{\alpha}: \alpha<\omega_{1}\}$ of elements of $H$ such that the sets of the family $\{g_{\alpha}+ Z : \alpha<\omega_{1}\}$ are pairwise disjoint. The following theorem generalizes Theorem K1 and Theorem K2 in which we assume Continuum hypothesis.\\
\vspace{.05cm}

\textbf{THEOREM 2.9 :} In every uncountable commutative Polish group $G$, there exist a proper $\sigma$-ideal $\mathcal I$ and a countable family $(B_{i})_{_{i<\omega}}$ of Bernstein subgroups such that if $\mathcal S$ is a $\sigma$-algebra and $\{\mathcal N_{n}\}_{_{n<\omega}}$ is a small system on $G$ satisfying the conditions:\\
(i) $\mathcal I\subseteq \mathcal N_{n}$ for all $n$\\
(ii) $\mathcal S$ is admissible with respect to $\{\mathcal N_{n}\}_{_{n<\omega}}$,\\
then at least one member of the above family of subgroups do not belong to the $\sigma$-algebra $\mathcal{\widetilde S}$ generated by $\mathcal S$ and $\mathcal N_{\infty}$.\\
\vspace{.05cm}

\textbf{PROOF :}  We first assume that $\{\mathcal N_{n}\}_{_{n<\omega}}$ is upper semicontinuous ralative to $\mathcal S$. We set $\mathcal T=\{ T\subseteq H : T+B\in \mathcal{\widetilde S}\}$ and\\
 $\mathcal M_{n}=\{ M\subseteq H : M+B\in \mathcal N_{n}\}$ $(n<\omega)$.\\
Then it can be easily checked that $\mathcal T$ is a $\sigma$-algebra and $\{\mathcal M_{n}\}_{_{n<\omega}}$ is small system of sets in $H$. Again as $x+B\subseteq \mathcal I\subseteq \mathcal N_{n}\subseteq \mathcal{\widetilde S}$ for all $n$ and $x\in H$, so $\mathcal T$ and $\{\mathcal M_{n}\}_{_{n<\omega}}$ are diffused.\
\vspace{.05cm}

By admissibility of $\mathcal S$ with respect to $\{\mathcal N_{n}\}_{_{n<\omega}}$, $\mathcal S\setminus \mathcal N_{n}\neq\emptyset$. Therefore $G\notin \mathcal N_{\infty}$ and consequently $H\notin \mathcal M_{\infty}$. Now $H$ being commutative, by Kulikov's theorem (see [8], [9]) we can write $H=\displaystyle\bigcup_{i<\omega}{\Gamma_{i}}$ where $\{\Gamma_{i}\}_{i<\omega}$ is an increasing sequence of commutative subgroups of $H$ and each $\Gamma_{i}$ is again a direct sum of cyclic groups. Thus for every $i<\omega$, $\Gamma_{i}= \displaystyle\sum_{j<\omega_{1}}[e_{ij}]$, where $[e_{i,j}]$ is the cyclic group generated by the element $e_{i,j}$ in $H$. Let $E_{i}=\{ e_{i,j} : j<\omega_{1}\}$ and $(E_{_{i,k,n}}){_{_{k<\omega_{1},n<\omega}}}$ be the Banach Kuratowski matrix over $E_{i}$. We put $\Gamma_{_{i,k,n}}= [E_{_{i,k,n}}]$ where for any set $E$, $[E]$ represents the group generated by $E$. Then $(\Gamma_{_{i,k,n}}){_{_{k<\omega,n<\omega}}}$ is a Banach Kuratowski matrix over $\Gamma_{i}$. We assert that the family $\{\Gamma_{i} : i<\omega\}\cup\{\Gamma_{_{i,k,n}} : i<\omega, k<\omega, n<\omega\}$ is not contained in $\mathcal T$.\
\vspace{.05cm}

Suppose if possible, let all the members of the above family belong to $\mathcal T$. Since $\mathcal M_{\infty}$ is a $\sigma$-ideal, then there exists $i^{\ast}$ such that $\Gamma_{i^{\ast}}\notin\mathcal M_{\infty}$.By Proposition 2.5, let $m$ be the natural number such that no subset $M$ of $\Gamma_{i^{\ast}}$ belongs to $\mathcal M_{m}$ if its complement in $\Gamma_{i^{\ast}}$ is in $\mathcal M_{m}$.\
\vspace{.05cm}

We set $G_{_{k,n}}= \Gamma_{i^{\ast}}\setminus \Gamma_{i^{\ast},k,n}$; $k,n<\omega$. Then $\displaystyle\bigcap_{n<\omega}{G_{k,n}}=\emptyset$ and therefore $\displaystyle\bigcap_{n<\omega}{(G_{k,n}+B)}=B\in\mathcal N_{\infty}$ where $\{G_{k,n}+B\}_{n<\omega}$ is a nested sequence of sets from $\mathcal{\widetilde S}$. Since $\{\mathcal N_{n}\}_{n<\omega}$ is upper semicontinuous relative to $\mathcal{\widetilde S}$ (Proposition 2.4), there exists $n_{k}$ such that $G_{k,n_{_{k}}}+ B\in\mathcal N_{n_{_{k}}}$. Hence $G_{n_{_{k}}}\in\mathcal M_{n_{_{k}}}$. By condition (iv) of Definition 2.1, we may choose $n_{k}>m$ such that $\displaystyle\bigcup_{k<\omega}{G_{n_{_{k}}}}\in\mathcal M_{m}$. Hence $\displaystyle\bigcap_{k<\omega}{\Gamma_{i^{\ast},k,n_{_{k}}}}\notin\mathcal M_{m}$. But $(\Gamma_{i^{\ast},k,n})_{k<\omega,n<\omega}$ being a Banach Kuratowski matrix, card$(\displaystyle\bigcap_{k<\omega}\Gamma_{i^{\ast},k,n_{_{k}}})\leq\omega$. So $\displaystyle\bigcap_{k<\omega}\Gamma_{i^{\ast},k,n_{_{k}}}\in\mathcal M_{\infty}$ because $\{\mathcal M_{n}\}_{n<\omega}$ is diffused. This is a contradiction. We denote the family $\{\Gamma_{i} : i<\omega\}\cup\{\Gamma_{i,k,n} : i<\omega, k<\omega, n<\omega\}$ by $\{G_{i} : i<\omega\}$. Consequently, not all members of the family $\{G_{i}+ B : i<\omega\}$ which consists of Bernstein groups (a fact which can be easily checked) belong to $\mathcal{\widetilde S}$. This proves the theorem.\
\vspace{.1cm}

The above proof rests heavily on the notion of upper semicontinuity of $\{\mathcal N_{n}\}_{n<\omega}$ relative to $\mathcal S$. If we remove this condition then we can no longer use Banach Kuratowski matrix. But still then we are able to prove the theorem by suitably constructing an admissible matrix on $H$. We show this below.\
\vspace{.1cm}

As before, we use Kulikov's theorem and write $H= \displaystyle\bigcup_{i<\omega}{\Gamma_{i}}$ where $\Gamma_{i}=\displaystyle\sum_{j<\omega_{1}}[e_{i,j}]$. Let $A_{i}= \{e_{i,j} : j<\omega_{1}\}$ and $(A_{i,n,\xi})_{_{n<\omega_{1},\xi<\omega_{1}}}$ be an admissible matrix over $A_{i}$. By lemma 6 (Chapter 13, [7]), for each $i$, there exists a countable family $\mathcal S_{i}$ of sets from $A_{i}$ such that $\{A_{i,n,\xi} : n<\omega_{1}, \xi<\omega_{1}\}\subseteq\sigma(\mathcal S_{i})$ where $\sigma(\mathcal S_{i})$ is the $\sigma$-algebra generated by $\mathcal S_{i}$. We may assume that $\mathcal S_{i}$ is an algebra, so that $\sigma(\mathcal S_{i})$ is the monotone class generated by $\mathcal S_{i}$. Then all the subgroups $[A_{i,n,\xi}]$ ($n<\omega, \xi<\omega_{1}$) belong to the monotone class generated by the family $\{[Z] : Z\in\mathcal S_{i}\}$. Since $\mathcal M_{\infty}$ is a $\sigma$-ideal, there exists $i^{\ast}$ such that $\Gamma_{i^{\ast}}\notin \mathcal M_{\infty}$. Moreover, $\{[A_{i^{\ast},n,\xi}] : n<\omega_{1}, \xi<\omega_{1}\}$ is an admissible matrix over $\Gamma_{i^{\ast}}$. Suppose all the members of the family $\{\Gamma_{i} : i<\omega\}\cup \{[Z] : Z\in\displaystyle\bigcup_{i<\omega}{\mathcal S_{i}}\}$ belong to $\mathcal T$. Then so are all the members of the family $\{[Z] : Z\in\mathcal S_{i^{\ast}}\}$ and consequently all the members $\{[A_{i^{\ast},n,\xi}] : n<\omega, \xi<\omega_{1}\}$ belong to $\mathcal T$. Hence there exists  a set $\Xi\subseteq\omega_{1}$ with card($\Xi$)=$\omega_{1}$ and $n_{_{0}}$ such that $[A_{i^{\ast},n_{_{0}},\xi}]\notin \mathcal M_{\infty}$ for all $\xi\in\Xi$ and $\bigcap\{[A_{i^{\ast},n_{_{0}},\xi}] : \xi\in\mathcal D\}= \{0\}$ for any set $\mathcal D\subseteq\omega_{1}$ having card($\mathcal D$)= $n_{_{0}}+2$. Consequently, $[A_{i^{\ast},n_{_{0}},\xi}]+ B\notin \mathcal N_{\infty}$ for all $\xi\in\Xi$ and $\bigcap\{[A_{i^{\ast},n_{_{0}},\xi}]+ B : \xi\in\mathcal D\}=\{0\}+B\in\mathcal N_{\infty}$. But this contradicts Proposition 2.7 because $\mathcal{\widetilde S}$ is admissible with respect to $\{\mathcal N_{n}\}_{n<\omega}$. Thus all the members of the above family cannot belong to $\mathcal T$. We denote this family by $\{G_{i} : i<\omega\}$. Consequently, not all members of the family $\{G_{i}+B : i<\omega\}$ which consists of Bernstein groups (a fact which can be easily checked) belong to $\mathcal{\widetilde S}$.
\begin{center}

\end{center}

\textbf{\underline{AUTHOR'S ADDRESS}}\

{\normalsize
\textbf{S.Basu}\\
\vspace{.02cm}
\hspace{.28cm}
\textbf{Dept of Mathematics}\\
\vspace{.02cm}
\hspace{.35cm}
\textbf{Bethune College, Kolkata} \\
\vspace{.01cm}
\hspace{.3cm}
\textbf{W.B. India}\\
\vspace{.02cm}
\hspace{.1cm}
\hspace{.29cm}\textbf{{e-mail : sanjibbasu08@gmail.com}}\\

\textbf{D.Sen}\\
\vspace{.02cm}
\hspace{.35cm}
\textbf{Saptagram Adarsha vidyapith (High), Habra, $\textbf{24}$ Parganas (North)} \\
\vspace{.01cm}
\hspace{.35cm}
\textbf{W.B. India}\\
\vspace{.02cm}
\hspace{.1cm}
\hspace{.29cm}\textbf{{e-mail : reachtodebasish@gmail.com}}
}

\end{document}